\documentclass[12pt]{amsart}

\usepackage{eucal,amsmath,amsthm,amsfonts,verbatim,amsopn,amssymb}
\usepackage{fullpage}

\usepackage[all]{xy}

\newcommand\p{\partial}
\newcommand\tint{{\textstyle\int}}

\newcommand\<{\langle}
\renewcommand\>{\rangle}
\newcommand\half{\tfrac{1}{2}}
\renewcommand\phi{\varphi}

\newcommand\lra{\longrightarrow}

\newcommand\R{{\mathbb{R}}}

\newcommand\Z{{\mathbb{Z}}}

\newcommand\CA{\mathcal A}

\newcommand\CD{\mathcal D} 
\newcommand\CE{\mathcal E}
\newcommand\CF{\mathcal F}

\newcommand\CH{\mathcal H}
\newcommand\CI{\mathcal I}

\newcommand\CL{\mathcal L}
 
\newcommand\CN{\mathcal N} 

\newcommand\CP{\mathcal P}
\newcommand\CQ{\mathcal Q}
\newcommand\CR{\mathcal R}
\newcommand\CS{\mathcal S}

\newcommand\CV{\mathcal V}

\newcommand{\HH}{{\mathbf{H}}}
\newcommand{\CAhat}{{\hat\CA}}
\newcommand{\CLhat}{{\hat{\CL}}}

\newtheorem{theorem}{Theorem}[section]

\newtheorem{lemma}[theorem]{Lemma}
\newtheorem{corollary}[theorem]{Corollary}
\newtheorem{proposition}[theorem]{Proposition}

\theoremstyle{definition}

\newcommand{\ignore}[1]{}

\renewcommand{\epsilon}{\varepsilon}

\title{On the moduli space of deformations of bihamiltonian
  hierarchies of hydrodynamic type}

\begin{document}
\author{Aliaa Barakat}

\address{Department of Mathematics\\
Massachusetts Institute of Technology}

\email{barakat@math.mit.edu}
\date{}

\maketitle

\begin{abstract}
We investigate the deformation theory of the simplest bihamiltonian
structure of hydrodynamic type, that of the dispersionless KdV hierarchy. 
We prove that all of its deformations are \emph{quasi-trivial} in the
sense of B.\ Dubrovin and Y.\ Zhang, that is, trivial after allowing
transformations where the first partial derivative $\p u$ of the field
is inverted.  We reformulate the question about deformations as a
question about the cohomology of a certain double complex, and
calculate the appropriate cohomology group.
\end{abstract}

\section{Introduction}

In the early 1990's, M.\ Kontsevich's solution \cite{K} of E.\ 
Witten's conjecture established an intriguing connection between the
differential equations of the Korteweg-de Vries (KdV) hierarchy and 
the quantum theory of
two-dimensional topological gravity.  Their work suggested a deep
relationship between integrable hierarchies ``of the KdV-type'' and a
wide class of quantum field theories (QFTs).

Integrability of a differential equation manifests itself in a set of
properties, one of which is the existence of a bihamiltonian
structure. This means that the system can be written in Hamiltonian
form in two distinct ways with respect to two compatible Poisson
brackets, permitting the construction of an infinite integrable
hierarchy of differential equations containing the original equation.
This in turn gives a method of generating infinitely many symmetries
and conservation laws.  Historically, the KdV equation was the first
differential equation found to have such a remarkable set of
properties. This equation appeared in the nineteenth century in the
study of solitary water waves, but since then has played a prominent
role in many distinct areas of mathematics.

In a series of pioneering papers, B.\ Dubrovin investigated relations
between the structure of two-dimensional topological QFTs (TFTs) and
the theory of bihamiltonian structures. Dubrovin showed that the tree
level, or genus zero, approximation of the TFT is described by a
bihamiltonian integrable hierarchy, the genus zero hierarchy, which
encodes all of the structure of the model in this approximation.

In this paper, we analyze bihamiltonian deformations of
integrable hierarchies and their bihamiltonian structures
that are related to generalizations of the
Witten conjecture.  Following the discoveries of Kontsevich and
Witten, Dubrovin and Y.\ Zhang \cite{D1,D2,D3,DZ1} conjectured the
existence of a bihamiltonian deformation of the genus zero hierarchy
which encodes the recursion relations among the correlators of the 2d
TFT. They further conjectured that this deformed hierarchy should be
quasi-trivial, i.e.\ obtained from the genus zero hierarchy by some
generalized coordinate change.

For example, the genus zero hierarchy associated to 2d topological
gravity with trivial background is the dispersionless KdV hierarchy.
It is known that the KdV hierarchy, which by Kontsevich's theorem
gives a complete description of 2d topological gravity with trivial
background, is a quasi-trivial bihamiltonian deformation of the
dispersionless KdV hierarchy.  The main theorem, which is presented in
Section 3, proves that \emph{all} bihamiltonian deformations of the
dispersionless KdV hierarchy are quasi-trivial.

Given a Poisson bracket
on functionals one can associate to it a Poisson bi-vector in the
Schouten Lie algebra of functional multi-vectors. This is a super Lie
algebra with respect to the Schouten bracket, and it becomes a
differential graded (dg) Lie algebra when endowed with the adjoint
action of the bi-vector.  The moduli space of deformations of the
Poisson bracket is then controlled by the cohomology of this dg Lie
algebra.

\smallskip
The goal of Section 2 is to develop a convenient setup for carrying
out cohomology computations.  Of central importance to our study is
the notion of the Schouten bracket for functional multi-vectors. 
The Schouten bracket is due to I.\ Dorfman and I.\ Gelfand,
but their definition is not well suited to calculations, especially
those involving higher multi-vectors. The Schouten graded Lie algebra
is isomorphic to the Gerstenhaber algebra of superfunctions on a
canonical symplectic supermanifold. This isomorphism identifies the
Schouten bracket with an odd Poisson bracket, and enables us to give
explicit formulas for the Schouten bracket. We also develop a normal
form for functional multi-vectors, which makes the study of
obstructions to deformations far more straightforward.  
We expect this normal form to have applications beyond the 
problems treated here.

Section 3 begins with a basic introduction to the deformation theory
of bihamiltonian structures, or equivalently, of a pair of compatible
Poisson brackets.  The study of the moduli space of deformations of a
bihamiltonian structure reduces to questions about the
\emph{bihamiltonian cohomology} of the pair.  This is the cohomology
of a double complex, introduced in Section~\ref{HamPair}.  Equivalent
classes of infinitesimal deformations are parametrized by the
cohomology classes of the first bihamiltonian cohomology group, and
the obstruction classes to continuing an infinitesimal deformation to
higher orders are elements of the second bihamiltonian cohomology
group.

In Section~\ref{InfSymm} we determine the space of infinitesimal
symmetries of the bihamiltonian structure of the dispersionless KdV
hierarchy.  In Section~\ref{QuasiTrivProof}, we prove that all
bihamiltonian deformations of this structure are quasi-trivial.
We obtain our results by calculating the corresponding 
cohomology groups. 

\smallskip
We follow the summation convention, taking sums over equal upper and lower 
indices.  This work constitutes part of the author's dissertation~\cite{B}. We have 
striven to give a self-contained  treatment.  Theorem~\ref{mainthm} has also 
been proven independently by Liu and Zhang \cite{LZ2}.

\subsection*{Acknowledgments}

I am grateful to my advisor Prof.\ Ezra Getzler for introducing me to
this beautiful subject. I take pleasure in acknowledging the AAUW
Foundation for its support in the academic year 2004-2005, during
which I held an AAUW dissertation fellowship.  During my graduate studies, 
I was also supported by a VIGRE grant at the University of Chicago, 
and NSF grant DMS-0072508.

\section{The Schouten bracket and Hamiltonian operators}

\subsection{Hamiltonian operators}\label{HamOp}

Let $(u^\alpha)$, $1\le \alpha\le q$, be coordinates on an open subset
$U$ of a $q$-dimensional manifold $M$. Let $\CA_0$ be the algebra of
smooth functions $C^\infty(U)$ in the variables $u^\alpha$, and let
$\CA$ be the algebra of differential polynomials in the jet variables
$u^\alpha_k$, $k\geq 1$, over the algebra $\CA_0$.

The algebra $\CA$ is equipped with a natural grading
$\CA=\bigoplus_k\CA_k$, defined by setting
\begin{equation}
  \label{grading}
  \deg u^\alpha_k = k .
\end{equation}
Elements of $\CA$ will be referred to as \emph{polynomial densities}.
Densities in $\CA_k$ are homogeneous of degree $k$.

Let $\CA[n]$ be the algebra of differential polynomials in the jet
variables $u^\alpha_k$, $k\le n$, over $\CA_0$. The sequence
\begin{equation*}
  0 \subset \CA[0] = \CA_0 \subset \CA[1] \subset \dots \subset \CA[n]
  \subset \dots \subset \CA
\end{equation*}
is an increasing filtration of $\CA$. If $g\in\CA[n]$, we say that $g$
has order $\le n$.

The derivation representing differentiation with respect to $x$ is
given by the \emph{total derivative}
\begin{equation*}
  \p=\sum_{j=0}^\infty u_{j+1}^\alpha\p_{u^\alpha_j}.
\end{equation*}
A derivation $X$ is \emph{evolutionary vector field} if it satisfies
the commutation relation
\begin{equation}
  \label{evol}
  [X,\p] = 0 .
\end{equation}
An evolutionary vector field has the form
\begin{equation*}
  X = \p^jf^\alpha \, \p_{u^\alpha_j} ,
\end{equation*}
where the $q$-tuple $(f^1,\dots,f^q)$ of polynomial densities,
$f^\alpha\in\CA$, is called the \emph{characteristic} of $X$.  It is
evident from the commutation relation~\eqref{evol} and the Jacobi
identity for the Lie bracket, that the evolutionary vector fields form
a Lie algebra, which we denote by $\CV^1$.

Since $X$ is uniquely determined by its characteristic, we can
identify the space $\CV^1$ of evolutionary vector fields with the
space $\CA^q$.

A \emph{differential operator} is an element $\CD\in\CA[\p]$; such an
operator defines an endomorphism of $\CA$. A differential operator
$\CD$ is \emph{homogeneous} of order $k$ if it has the form
\begin{equation*}
  \CD = \sum_{j=0}^k P_j \p^j, \quad P_j\in\CA_{k-j} .
\end{equation*}
The \emph{adjoint} $\CD^*$ of a differential operator $\CD$ is given
by the formula
\begin{equation*}
  \CD^* = \sum_{j=0}^k \sum_{i=0}^{k-j} (-1)^{j+i} \tbinom{j+i}{i}
  (\p^iP_{j+i}) \p^j .
\end{equation*}

Define the space of \emph{functionals} by
\begin{equation*}
  \CF = \CA/\p\CA .
\end{equation*}
The \emph{variational derivative}
\begin{equation*}
  \delta_{u^\alpha} : \CF \to \CA
\end{equation*}
is given by the formula
\begin{equation}\label{EulerOp}
  \delta_{u^\alpha} = \sum_{k=0}^\infty (-\p)^k \p_{u^\alpha_k} .
\end{equation}

A $(q\times q)$-matrix $\CD^{\alpha\beta}$ of differential operators
defines a bracket
\begin{equation*}
  \{F,G\}_\CD = \sum_{\alpha,\beta} (\delta_{u^\alpha}F)
  \CD^{\alpha\beta} (\delta_{u^\beta}G) \quad \textup{mod $\p\CA$}
\end{equation*}
on functionals. This bracket is skew-symmetric if $\CD$ is
skew-adjoint, that is,
\begin{equation*}
  (\CD^{\alpha\beta})^* = - \CD^{\beta\alpha} .
\end{equation*}
If in addition, the bracket satisfies the Jacobi identity, then $\CD$
is called a \emph{Hamiltonian operator} or a \emph{Hamiltonian
  structure}. A pair of Hamiltonian operators $\CP$ and $\CQ$ is said
to form a \emph{bihamiltonian structure}, if for any scalar $\lambda$
the sum
\begin{equation*}
  \CP + \lambda\CQ
\end{equation*} 
is again a Hamiltonian operator.  The operators $\CP$ and $\CQ$ are
also called \emph{compatible}, and the brackets $\{\cdot,\cdot\}_\CP$
and $\{\cdot,\cdot\}_\CQ$ are said to form a \emph{Poisson pencil}.
 
The Schouten graded Lie algebra, introduced by Gelfand and Dorfman
\cite{GD}, is a graded vector space $\CV=\bigoplus_k\CV^k$ with super
Lie bracket
\begin{equation*}
  [[\cdot,\cdot]]:\CV^{k_1}\times\CV^{k_2}\lra\CV^{k_1+k_2-1} ,
\end{equation*}
such that $\CV^0\cong\CF^0$, $\CV^1$ is the space of evolutionary
vector fields, and $\CV^2$ is isomorphic to the space of skew-adjoint
$(q\times q)$ matrices of differential operators.  We recall the
precise definition of the bracket $[[\cdot,\cdot]]$ in the next
section. We conclude this section with the following standard
characterizations of Hamiltonian and bihamiltonian structures.

\begin{proposition}\label{HamCond}
  For a bi-vector $\CD\in\CV^2$, the following conditions are
  equivalent:
  \begin{enumerate}
  \item $\CD$ is Hamiltonian;
  \item the associated bracket $\{\cdot,\cdot\}_\CD$ is Poisson;
  \item the morphism $d_\CD=[[\CD,\cdot]]$ is a differential on $\CV$.
  \end{enumerate}
\end{proposition}

\begin{proposition}
  Two Hamiltonian operators $\CP$ and $\CQ$ are compatible if and only if 
  \begin{equation*}
    [[\CP,\CQ]]=0.
  \end{equation*}
\end{proposition}

\subsection{The Schouten bracket} \label{NSbra}
In this section we give a formula for the Schouten bracket due to Getzler 
\cite{G}.  Let $T^*[1]M$ be the $\Z$-graded manifold underlying the graded vector
bundle over $M$ whose fiber at $u\in M$ is the cotangent space
$T^*_uM$ concentrated in degree $-1$. Let $(\theta_\alpha)$ be the
coordinates along the fibers of $T^*[1]U$ dual to the coordinates
$(u^\alpha)$ on $U$.

Let $\Lambda^\bullet_\infty=\CA[\theta_{\alpha,k}\mid k\geq 0]$ be the
exterior algebra of functions over $\CA$ with generators
$\theta_{\alpha,k}$ in degree $1$. Denote by $\p_{\theta_{\alpha,k}}$
the graded derivation (partial derivative) such that
\begin{equation*}
  \p_{\theta_{\alpha,k}} \theta_{\beta,\ell} = \delta^\alpha_\beta
  \delta_{k,\ell} .
\end{equation*}
Extend the total derivative $\p$ on $\CA$ to the algebra
$\Lambda^\bullet_\infty$ by the formula
\begin{equation*}
  \p = \sum_{k=0}^\infty \bigl( u^\alpha_{k+1} \p_{u^\alpha_k} +
    \theta_{\alpha,k+1} \p_{\theta_{\alpha,k}} \bigr) .
\end{equation*}
Denote by $\delta_{\theta_\alpha}$ the variational derivative on
$\Lambda_\infty^\bullet$, defined by the formula analogous to
\eqref{EulerOp},
\begin{equation*}
  \delta_{\theta_\alpha} = \sum_{k=0}^\infty (-\p)^k \p_{\theta_{\alpha,k}} .
\end{equation*}

Let
\begin{equation*}
  \CV^\bullet = \Lambda^\bullet_\infty/\p\Lambda^\bullet_\infty .
\end{equation*}
Elements of $\CV^k$ are denoted $\int F\,dx$, $F\in\Lambda_\infty^k$,
and are called \emph{functional $k$-vectors}, or just $k$-vectors, for
short. Notice that, by definition, $\CV^0=\CF$, where $\CF$ is the
space of functionals introduced in the previous section.

The graded vector space $\CV^\bullet$ is a graded Lie algebra, with
respect to the \emph{Schouten bracket}
\begin{equation*}
  [[\cdot,\cdot]]: \CV^{k_1} \times \CV^{k_2} \lra \CV^{k_1+k_2-1}
\end{equation*}
defined by the formula
\begin{align*}
  [[\tint F\,dx,\tint G\,dx]] &= \sum_{k,\ell} \tint \Big(
  (-1)^{|F|+1}\p^k \p_{\theta_{\alpha,\ell}}F \cdot \p^\ell
  \p_{u^\alpha_k} G -\p^k\p_{u^\alpha_\ell}F \cdot
  \p^\ell\p_{\theta_{\alpha,k}}G \Big) \,dx \\
  &= \tint \Big( (-1)^{|F|+1} \delta_{\theta,\alpha} F \cdot
  \delta_{u^\alpha} G - \delta_{u^\alpha}F \cdot
  \delta_{\theta_\alpha} G \Bigr) \, dx .
\end{align*}
For a proof that this bracket satisfies the graded Jacobi identity, see \cite{G}.

In order to have a graded Lie algebra whose bracket has degree zero, it is
convenient to replace the graded vector space $\CV^\bullet$ by its
shifted version
\begin{equation*}
  \CL^\bullet = \CV^{\bullet+1} .
\end{equation*}
A Hamiltonian operator $\CH$ is then an element of $\CL^1$ satisfying
the Maurer-Cartan equation $[[\CH,\CH]]=0$. For such an $\CH$, the
morphism $d_{\CH}=[[\CH,\cdot]]$ on $\CL$ is a derivation of degree
$1$, and $(\CL,d_{\CH})$ is a dg Lie algebra.

The definition of degree of homogeneity on $\CA$ extends to
$\Lambda^\bullet_\infty$, by setting
\begin{equation*}
  \deg \theta_{\alpha,k} = k , \quad k \ge 0.
\end{equation*}
Denote by $\CL^k\<\ell\>$ the space of $(k+1)$-vectors homogeneous of
degree $k+\ell$.

\subsection{The graded Lie algebra $\CLhat$}

In this section, take $q=1$, so there is only a single dependent
variable $u$. In the study of deformations of the dispersionless KdV
hierarchy, an important role is played by the algebra $\CAhat =
\CA[u_1^{-1}]$. The grading of the algebra $\CA$ extends to a grading
of $\CAhat$ if we assign to $u_1^{-1}$ homogeneity degree $-1$:
\begin{equation*}
  \CAhat = \bigoplus_k \CAhat_k
\end{equation*}
If $n>0$, let $\CAhat[n]=(\CA[n])[u_1^{-1}]$.

In parallel to the definition of $\CAhat$, let
$\hat\Lambda^\bullet_\infty=\Lambda^\bullet_\infty[u_1^{-1}]$. The
derivation $\p$ maps $\hat\Lambda^\bullet_\infty$ to itself, and we
define the graded Lie algebra $\CLhat^\bullet$ by
\begin{equation*}
  \CLhat^k = \hat\Lambda^{k+1}_\infty/\p\hat\Lambda^{k+1}_\infty .
\end{equation*}
The bracket $[[\cdot,\cdot]]$ on $\CLhat$, defined by the same
formulas as for the Schouten bracket on $\CL$, turns it into a graded
Lie algebra. Denote by $\CLhat^k\<\ell\>$ the space of $(k+1)$-vectors
homogeneous of degree $k+\ell$. The Schouten bracket is homogeneous of
degree $0$, in the sense that
\begin{equation*}
  [[ \CLhat^{k_1}\<\ell_1\> , \CLhat^{k_2}\<\ell_2\> ]] \subset
  \CLhat^{k_1+k_2}\<\ell_1+\ell_2\> .
\end{equation*}

\newcommand{\PP}{{\mathbf{P}}}
\newcommand{\QQ}{{\mathbf{Q}}}

\subsection{Formal deformations of Hamiltonian operators}

Let ${\CH}\in{\CL}^1$ be a Hamiltonian operator. 
Recall that $\CH$ being Hamiltonian is equivalent to the vanishing
of the Schouten bracket $[[{\CH},{\CH}]]$.

By a \emph{formal deformation} of $\CH$ we mean a formal power series
\begin{equation*}
  \HH={\CH}+\sum_{k=1}^\infty\epsilon^k{\CH}_k, \quad
  {\CH}_k\in{\CL}^1
\end{equation*}
such that $[[\HH,\HH]]=0$. 
The $n$-th order
deformations of $\CH$ are given by the Maurer-Cartan elements of the
dg Lie algebra ${\CL}\otimes I_n$, where $I_n$ is the ideal in
$\R[\epsilon]/(\epsilon^{n+1})$ generated by $\epsilon$.  They are of
the form
\begin{equation*}
  \HH_n={\CH}+\sum_{k=1}^n\epsilon^k{\CH}_k, \quad
  {\CH}_k\in{\CL}^1  
\end{equation*}
and the associated brackets on functionals satisfy the Jacobi identity up to order $n$.

Every first order deformation of ${\CH}$ determines, and is uniquely
determined, by an \emph{infinitesimal deformation}, i.e.\ a bi-vector
${\CH}_1\in{\CL}^1$ such that
\begin{equation}\label{InfDef}
  [[{\CH},{\CH}_1]]=0.
\end{equation}

The Schouten dg Lie algebra $({\CL},d_{\CH})$ controls the moduli
space of deformations of $\CH$ in the following sense.  A formal
deformation $\HH=\sum_{k=0}^\infty\epsilon^k{\CH}_k$ is called
trivial if there is a formal coordinate change
\begin{equation*}
  u^\beta\mapsto\psi^\beta [u^\alpha]=u^\beta+\epsilon f^\beta+\ldots  
\end{equation*}
such that
\begin{equation*}
  \HH={\CH} - \epsilon [[X,{\CH}]]+\ldots,
\end{equation*}
where $X=\tint f^\alpha \theta_\alpha \, dx$ is the evolutionary
vector field with characteristic $(f^\alpha)$.  In particular,
\begin{equation*}
  {\CH}_1=[[{\CH},X]]=d_{\CH}X.
\end{equation*}

If ${\CH}_1$ is an infinitesimal deformation, then the first
\emph{obstruction cocycle} is given by
\begin{equation*}
  [[{\CH}_1,{\CH}_1]]\in {\CL}^2.  
\end{equation*}
It is the coefficient of $\epsilon^2$ in $[[\HH_1,\HH_1]]$, and it
satisfies
\begin{align*}
  d_{\CH}\left([[{\CH}_1,{\CH}_1]]\right) &=
  [[{\CH},[[{\CH}_1,{\CH}_1]]\,]] \\
  &=
  [[\,\underbrace{[[{\CH},{\CH}_1]]}_{=0},{\CH}_1]]-[[{\CH}_1,\underbrace{[[{\CH},{\CH}_1]]}_{=0}\,]] \\
  &= 0 .
\end{align*}
Further, if $d_{\CH}{\CH}_2=-\half[[{\CH}_1,{\CH}_1]]$,
for some bi-vector ${\CH}_2\in{\CL}^1$, then
\begin{equation*}
  \HH_2={\CH}+\epsilon{\CH}_1+\epsilon^2{\CH}_2 
\end{equation*}
is a second order deformation.  For an $n$-th order deformation
$\HH_n$, the $n$-th obstruction cocycle is the coefficient of
$\epsilon^{n+1}$ in $[[\HH_n,\HH_n]]$. It is given by
\begin{equation*}
  \sum_{i=1}^n[[{\CH}_i,{\CH}_{n-i+1}]],  
\end{equation*}
and by a calculation parallel to the one above, we see that $\HH_n$
can be extended to an $(n+1)$-st order deformation
$\HH_{n+1}=\HH_n+\epsilon^{n+1}{\CH}_{n+1}$ if
\begin{equation*} 
  d_{\CH}{\CH}_{n+1}=-\half
  \sum_{i=1}^n[[{\CH}_{i},{\CH}_{n-i+1}]] .
\end{equation*}

To summarize, we have thus identified the equivalent infinitesimal
deformations of the Hamiltonian operator $\CH$ with the elements of
the first cohomology group of the dg Lie algebra $({\CL},d_{\CH})$,
and the obstruction classes to extending a first order deformation of
$\CH$ to higher order ones with elements of the second cohomology
group.

\subsection{Hamiltonian operators of hydrodynamic type}

A Hamiltonian operator $\CH$ has \emph{hydrodynamic type} if it is
homogeneous of order $1$ (Dubrovin and Novikov \cite{DN}):
\begin{equation*}
  \CH^{\alpha\beta} = h^{\alpha\beta}(u) \p
  + \Gamma^{\alpha\beta}_\gamma(u) u^\gamma_1 .
\end{equation*}
Skew-symmetry implies that $h^{\alpha\beta}$ is symmetric, and that
\begin{equation*}
  \Gamma^{\alpha\beta}_\gamma(u) + \Gamma^{\beta\alpha}_\gamma(u) =
  \p_\gamma h^{\alpha\beta}(u) .
\end{equation*}
Under a change of coordinates, the coefficients $h^{\alpha\beta}$
transform as the components of a symmetric bilinear form on the
cotangent bundle, and the coefficients
$\Gamma^{\alpha\beta}_\gamma(u)$ transform as the associated
Christoffel symbols
\begin{equation*}
  \Gamma^{\alpha\beta}_\gamma = -\half h^{\alpha\delta}
  h^{\beta\epsilon} \bigl( \p_\delta h_{\epsilon\gamma} - \p_\epsilon
  h_{\delta\gamma} + \p_\gamma h_{\delta\epsilon} \bigr) .
\end{equation*}
The condition that $\CH$ is a Hamiltonian operator may be expressed as
the flatness condition
\begin{equation*}
  h^{\alpha\delta} \p_\delta \Gamma^\beta - h^{\beta\delta} \p_\delta
  \Gamma^\alpha + [ \Gamma^\alpha , \Gamma^\beta ] = 0 .
\end{equation*}

We say that the Hamiltonian operator $\CH$ is \emph{nondegenerate} if
the associated bilinear form $h^{\alpha\beta}$ is nondegenerate, that
is, a flat pseudo-metric. In this case, we may choose flat
coordinates, i.e.\ those for which the components $h^{\alpha\beta}$ of
the pseudo-metric are constant. In terms of these, the Hamiltonian
operator $\CH$ becomes
\begin{equation*}
  \CH^{\alpha\beta} = h^{\alpha\beta} \p .
\end{equation*}

If $\CH$ is a Hamiltonian operator of hydrodynamic type, then the
differential $d_{\CH}$ on $\CL$ raises the degree by $1$, and hence
maps $\CL^k\<\ell\>$ to $\CL^{k+1}\<\ell\>$. Thus the dg Lie algebra
$(\CL,d_{\CH})$ decomposes into the direct sum of subcomplexes
\begin{equation*}
  \CL\<\ell\> = \bigoplus_k \CL^k\<\ell\> 
\end{equation*}
and the cohomology of $(\CL,d_{\CH})$ decomposes accordingly into
\begin{equation*}
  H^k(\CL,d_{\CH}) = \bigoplus_{\ell\ge-k} H^k(\CL\<\ell\>,d_{\CH}) .
\end{equation*}

The next theorem is due to Getzler \cite{G}.
\begin{theorem}
  \label{Getzler}
  Let $\CH$ be a Hamiltonian operator of hydrodynamic type. If $\CH$
  is nondegenerate,
  \begin{equation*}
    H^k(\CL\<\ell\>,d_{\CH})=0, \quad \text{if $\ell>-k$.}
  \end{equation*}
\end{theorem}

Now suppose that $\CH$ is a nondegenerate Hamiltonian operator of
hydrodynamic type, and let $\HH$ be a formal deformation of $\CH$,
\begin{equation*}
  \HH = {\CH} + \sum_{k=1}^\infty \epsilon^k {\CH}_k ,
\end{equation*}
such that the infinitesimal deformation ${\CH}_1$ is homogeneous of
degree $p+1$, that is, $\CH_1\in\CL^1\<p\>$, and for $k>1$,
\begin{equation*}
  \CH_k \in \bigoplus_{\ell\ge0} \CL^1\<\ell\> .
\end{equation*}

We may show, by induction on $k$, that the deformation $\HH$ is
equivalent to a deformation such that $\CH_k$ is homogeneous of degree
$kp+1$, that is, $\CH_k\in\CL^1\<kp\>$. To see this, observe that the
right-hand side of the equation
\begin{equation*}
  d_{\CH}{\CH}_k = - \half \sum_{i=1}^{k-1} [[{\CH}_i,{\CH}_{k-i}]] ,
\end{equation*}
lies in $\CL^2\<kp\>$. By Theorem \ref{Getzler}, there is an element
$\CI_k$ of $\CL^0$ such that $\CH_k+d_\CH\CI_k$ is homogeneous of
degree $kp+1$.

A deformation $\HH$ with this property is called homogeneous.  In this
paper, we only consider homogeneous deformations. 
The second Hamiltonian structure 
\begin{equation*}
  u\p + \half u_1 + \tfrac{3}{2} \epsilon^2 \p^3  
\end{equation*}
of the small dispersion expansion
\begin{equation}
  \label{kdv}
  \p_tu = uu_1 + \epsilon^2 u_3
\end{equation}
of the KdV equation is an example of such a deformation.

\subsection{The normalization operator} \label{normalform}

The \emph{higher variational derivatives} on $\Lambda^\bullet_\infty$
are defined by
\begin{equation*}
  \delta_{k,u^\alpha} = \sum_{j=0}^\infty (-1)^j \tbinom{k+j}{k}
  \p^j\p_{u^\alpha_{k+j}} \quad \text{and} \quad
  \delta_{k,\theta_\alpha} = \sum_{j=0}^\infty (-1)^j \tbinom{k+j}{k}
  \p^j\p_{{\theta_{\alpha,k+j}}} .
\end{equation*}
Note that $\delta_{0,u^\alpha}=\delta_{u^\alpha}$ and
$\delta_{0,\theta_\alpha}=\delta_{\theta_\alpha}$ are the variational
derivatives introduced earlier.

Let $\CN:\Lambda_\infty^\bullet \to \Lambda_\infty^\bullet$ be the
\emph{normalization operator}
\begin{equation*}
  \CN = \sum_\alpha \theta_\alpha \delta_{\theta_\alpha} .
\end{equation*}
Note that $\CN\cdot\p=0$, since $\delta_{\theta_\alpha}\p=0$.

\begin{theorem}
  If $F\in\Lambda^k_\infty$, $\CN F - k\,F \in \p \Lambda^k_\infty$.
\end{theorem}
\begin{proof}
  A direct calculation shows that
  \begin{equation*}
    \p_{\theta_{\alpha,i}} = \sum_{j=i}^\infty \tbinom{j}{i}
    \p^{j-i} \delta_{j,\theta_\alpha} .
  \end{equation*}
  It follows that if $F\in\Lambda^k_\infty$,
  \begin{align*}
    k\,F & = \sum_{j=0}^\infty \theta_{\alpha,j}\p_{\theta_{\alpha,j}} F \\
    & = \sum_{j=0}^\infty \sum_{i=0}^j \tbinom{j}{i} \theta_{\alpha,i}
    \p^{j-i}\delta_{j,\theta_\alpha}F \\
    & = \sum_{j=0}^\infty \p^j \bigl( \theta_\alpha
    \delta_{j,\theta_\alpha} F \bigr) \\
    & = \CN F + \p \cdot \sum_{j=0}^\infty \p^j \bigl( \theta_\alpha
    \delta_{j+1,\theta_\alpha} F \bigr) .
    \qedhere
  \end{align*}
\end{proof}

We will also use the generalization of this theorem with
$\Lambda^\bullet_\infty$ replaced by $\hat\Lambda^\bullet_\infty$; the
proof is identical.

\smallskip
The normalization operator $\CN$ may be used to establish the standard
identifications of the spaces $\CL^0$ and $\CL^1$ mentioned in
Section~\ref{NSbra}. For example, if
\begin{equation*}
  F = \sum_j f^\alpha_j \theta_{\alpha,j} \in \Lambda^1_\infty , \quad
  f^\alpha_j \in \CA ,
\end{equation*}
we have $\tint F\, dx = \tint \CN F\, dx \in \CL^0$, where the
normalization $\CN F$ is given by the formula
\begin{align*}
  \CN F &= \sum_j \theta_\alpha \delta_{\theta_\alpha} ( f^\alpha_j
  \theta_{\alpha,j} ) \\
  &= \sum_j \theta_\alpha (-\p)^j f^\alpha_j .
\end{align*}
That is, $\tint F\,dx$ is the evolutionary vector field with
characteristic $\bigl( \sum_j (-\p)^j f^\alpha_j \bigr)$.

Similarly, if
\begin{equation*}
  G = \sum_{j,k} f^{\alpha\beta}_{jk} \theta_{\alpha,j}
  \theta_{\beta,k} \in \Lambda^2_\infty ,
\end{equation*}
where $f^{\alpha\beta}_{jk}=-f^{\beta\alpha}_{kj}$, we have $\tint G
\, dx = \half \tint \CN G \, dx \in \CL^1$. Here, the normalization
$\CN G$ is given by the formula
\begin{equation*}
  \half \CN G = \sum_{j=0}^\infty \theta_\alpha (-\p)^j \bigl(
  f^{\alpha\beta}_{jk} \theta_{\beta,k} \bigr)
  = \theta_\alpha\CD^{\alpha\beta}\theta_\beta ,
\end{equation*}
where $\CD=(\CD^{\alpha\beta})$ is a skew-adjoint $(q\times q)$-matrix
of differential operators
\begin{equation*}
  \CD^{\alpha\beta} = \sum_{j=0}^\infty \sum_{i=0}^j (-1)^j
  \tbinom{j}{i} \bigl( \p^{j-i} f^{\alpha\beta}_{jk} \bigr) \p^{k+i} .
\end{equation*}

The role of the normalization operator $\CN$ is to implement analogues
of these identifications of $\CL^k$ for arbitrary $k\ge0$.

\section{The Bihamiltonian Cohomology for Hamiltonian Structures of
  Hydrodynamic Type}

\subsection{Bihamiltonian structures of hydrodynamic type}

The dispersionless KdV equation
\begin{equation}
  \label{disp}
  \p_t u = uu_1
\end{equation}
is contained in a sequence of commuting flows of one dependent variable
$u$ which make up the dispersionless KdV hierarchy. This hierarchy is
associated to the following compatible pair of Hamiltonian operators
of hydrodynamic type:
\begin{equation}
  \label{biham}
  \CP = \p \quad \text{and} \quad \CQ = u\p + \half u_1 .
\end{equation}
For the choice of Hamiltonian functionals $H_0= \tfrac{1}{3} \int u^2\, dx$ and
$H_1= \tfrac{1}{6} \int u^3\, dx$, this equation may be written
\begin{equation*}
  \p_t u = \CP\delta_u H_1 = \CQ\delta_u H_0.  
\end{equation*}
The pseudo-differential operator
\begin{align*}
  \CR &= \CQ \CP^{-1} \\
  &= u +\half u_1\p^{-1}
\end{align*}
is a recursion operator \cite{O} for the dispersionless KdV equation.
Letting
\begin{equation*}
  \p_{t_0} u = \CP \delta_u H_0 \quad \text{and} \quad \p_{t_1} u =
  \p_t u = \CP \delta_u H_1 ,
\end{equation*}
the equations of the dispersionless KdV hierarchy are given
recursively by
\begin{equation*}
  \p_{t_n} u = \CP\delta_u H_n = \CQ\delta_u H_{n-1}
\end{equation*}
for Hamiltonian functionals $H_n$, $n\geq 0$.  The functionals can be
recursively determined by the relations
$\{\cdot,H_n\}_{\CP}=\{\cdot,H_{n-1}\}_{\CQ}$ starting with the
Casimir $H_{-1}=\tfrac{4}{3} \int u \, dx$ for $\{\cdot,\cdot\}_{\CP}$.  They are
in involution with respect to both Poisson brackets
$\{\cdot,\cdot\}_{\CP}$ and $\{\cdot,\cdot\}_{\CQ}$.

We can now formulate the problem posed by Dubrovin and Zhang
\cite{DZ2}. Suppose we are given a compatible pair $(\CP,\CQ)$ of
Hamiltonian operators of hydrodynamic type, such that $\CP$ is
nondegenerate. One wishes to classify homogeneous deformations of this
pair, modulo the action of the group of \emph{Miura transformations},
that is, coordinate changes of the form
\begin{equation*}
  u^\alpha \mapsto u^\alpha + \sum_{k=1}^\infty \epsilon^k
  \psi^\alpha_k , \quad \psi^\alpha_k \in \CA_k .
\end{equation*}

In this work, we show that every non-trivial deformation of the
bihamiltonian structure~\eqref{biham} of the dispersionless KdV
hierarchy can be transformed into $(\CP,\CQ)$ by a \emph{quasi-Miura
  transformation}, that is, a coordinate change of the form
\begin{equation}
  \label{qmiura}
  u \mapsto u + \sum_{k=1}^\infty \epsilon^k \psi_k , \quad
  \psi_k \in \CAhat .
\end{equation}
We say that every homogeneous deformation of $(\CP,\CQ)$ is
\emph{quasi-trivial}.

Our work was motivated by Lorenzoni \cite{L}: he studied homogeneous
deformations of this bihamiltonian structure, and showed by explicit
calculation that they were quasi-trivial up to fourth order in $\epsilon$.

We conclude this section by showing that it is indeed necessary to
consider quasi-Miura transformations. Suppose that $u$ satisfies the
dispersionless KdV equation~\eqref{disp}, and that
\begin{equation*}
  v = u + \epsilon^2 \psi + O(\epsilon^4)
\end{equation*}
satisfies the KdV equation~\eqref{kdv} up to terms of order
$\epsilon^4$. We calculate that
\begin{align*}
  \p_tv - v\p v - \epsilon^2 \p^3v &= \bigl( \p_tu + \epsilon^2
  \p_t\psi \bigr) - \bigl( uu_1 + \epsilon^2 u\p\psi + \epsilon^2
  u_1\psi \bigr) - \epsilon^2 u_3 + O(\epsilon^4) \\
  &= \epsilon^2 \bigl( \p_t\psi - u\p\psi - u_1\psi - u_3 \bigr) +
  O(\epsilon^4) = 0 .
\end{align*}
This equation for $\psi$ has no solution in $\CA$, but it does have a
solution in $\CAhat$, namely
\begin{align*}
  \psi &= - \frac{1}{2} \Bigl( \frac{u_3}{u_1} - \frac{u_2^2}{u_1^2}
  \Bigr) \\
  &= - \half \, \p^2 \log(u_1) \in \CAhat[3] .
\end{align*}

\subsection{Formal deformations of a compatible pair of Hamiltonian
  operators} \label{HamPair}

Let $\CP$ and $\CQ$ be two compatible Hamiltonian operators. We
introduce a double complex
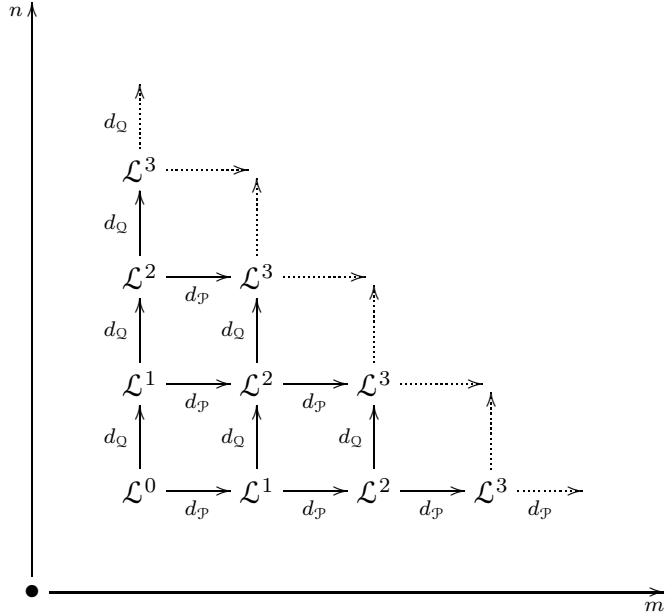
\begin{equation*}
  C^{\bullet,\bullet} = \bigoplus C^{m,n}
\end{equation*}
\begin{equation*}
  \xymatrix{
    & & & & & & & \\
    & & & & & & & \\
    & {\CL}^3 \ar@{.>}[u]^{d_{\CQ}} \ar@{.>}[r] & & & & & & \\
    & {\CL}^2 \ar[u]^{d_{\CQ}} \ar[r]_{d_{\CP}} & {\CL}^3 \ar@{.>}[u]
    \ar@{.>}[r] & & & & & \\
    & {\CL}^1 \ar[u]^{d_{\CQ}} \ar[r]_{d_{\CP}} & {\CL}^2
    \ar[u]^{d_{\CQ}} \ar[r]_{d_{\CP}} & {\CL}^3 \ar@{.>}[u] \ar@{.>}[r]
    & & & & \\
    & {\CL}^0 \ar[u]^{d_{\CQ}} \ar[r]_{d_{\CP}} & {\CL}^1
    \ar[u]^{d_{\CQ}} \ar[r]_{d_{\CP}} & {\CL}^2 \ar[u]^{d_{\CQ}}
    \ar[r]_{d_{\CP}} & {\CL}^3 \ar@{.>}[r]_{d_{\CP}} \ar@{.>}[u] & & & \\
    \bullet \ar[rrrrrr]_>>{m}\ar[uuuuuu]^>>{n} & & & & & & &
  }
\end{equation*}
with $C^{m,n}={\CL}^{m+n}$. To see that this is indeed a double
complex, observe that the differentials $d_{\CP}$ and $d_{\CQ}$
anticommute: if $\CR$ is an element of $\CL$,
\begin{align*}
  d_{\CP}d_{\CQ}(\CR)&= [[{\CP},[[{\CQ},{\CR}]]\,]] \\
  &=
  [[\,\underbrace{[[{\CP},{\CQ}]]}_{=0},{\CR}]]-[[{\CQ},[[{\CP},{\CR}]]\,]] \\
  &= -d_{\CQ}d_{\CP}({\CR}) .
\end{align*}
The associated total complex $C^\bullet$ with differential
$d=d_{\CP}+d_{\CQ}$ is obtained by summing along the anti-diagonals:
\begin{equation*}
  C^k = \bigoplus_{m+n=k} C^{m,n} = \bigoplus_{m+n=k} \CL^{m+n} .
\end{equation*}
The \emph{bihamiltonian cohomology} 
\begin{equation*}
  H^\bullet({\CL};d_{\CP},d_{\CQ}) =
  Z^\bullet({\CL};d_{\CP},d_{\CQ})/B^\bullet({\CL};d_{\CP},d_{\CQ})
\end{equation*}
is the cohomology of this total complex.

A cocycle in $Z^k({\CL};d_{\CP},d_{\CQ})$ is a $(k+1)$-tuple
\begin{equation*}
  (c_0,\ldots,c_k) \in \CL^k \oplus \dots \oplus \CL^k ,
\end{equation*}
such that $d_{\CP}c_0=0$, $d_{\CQ}c_i+d_{\CP}c_{i+1}=0$ for $0\le i<k$,
and $d_{\CQ}c_k=0$. This cocycle is a coboundary if there is a
$(k-1)$-cochain $(a_0,\ldots,a_{k-1})$ such that $c_0=d_{\CP}a_0$,
$c_i=d_{\CQ}a_{i-1}+d_{\CP}a_i$ for $0<i<k$, and $c_k=d_{\CQ}a_{k-1}$.

Now suppose that $\PP={\CP}+\sum_{k=1}^\infty\epsilon^k{\CP}_k$
and $\QQ={\CQ}+\sum_{k=1}^\infty\epsilon^k{\CQ}_k$ are two
\emph{compatible} formal deformations of $\CP$ and $\CQ$.  Writing the
differential $d=d_{\CP}+d_{\CQ}$ in matrix form, we calculate that
\begin{align*}
  d({\CP}_1,{\CQ}_1) &=
  \begin{pmatrix}
    d_{\CP} & 0 \\
    d_{\CQ} & d_{\CP} \\
    0 & d_{\CQ}
  \end{pmatrix}
  \begin{pmatrix}
    {\CP}_1 \\
    {\CQ}_1
  \end{pmatrix} \\
  &=
  \begin{pmatrix}
    d_{\CP}{\CP}_1 \\
    d_{\CQ}{\CP}_1+d_{\CP}{\CQ}_1 \\
    d_{\CQ}{\CQ}_1
  \end{pmatrix}
  \in {\CL}^2 \oplus \CL^2 \oplus \CL^2 .
\end{align*}
The entries of the latter vector are the coefficients of $\epsilon$ in
the compatibility conditions $\half[[\PP,\PP]]=0$, $[[\PP,\QQ]]=0$ and
$\half[[\QQ,\QQ]]=0$, respectively.

A first order bihamiltonian deformation consists of first order
deformations
\begin{equation*}
  \PP_1={\CP}+\epsilon{\CP}_1 \qquad \text{and} \qquad
  \QQ_1={\CQ}+\epsilon{\CQ}_1
\end{equation*}
compatible up to first order. The first obstruction cocycle is given
by
\begin{equation*}
  \bigl( \half[[{\CP}_1,{\CP}_1]] ,
  [[{\CP}_1,{\CQ}_1]] ,
  \half[[{\CQ}_1,{\CQ}_1]]) \in \CL^2 \oplus \CL^2 \oplus
  \CL^2 ,
\end{equation*}
and one easily sees that it vanishes under the differential $d:C^2\to
C^3$. For example, using that
$[[{\CP},{\CQ}_1]]+[[{\CQ},{\CP}_1]]=0$, the graded Jacobi
identity for the Schouten bracket, and the graded anti-commutativity,
we calculate that
\begin{align*}
  d_{\CP}[[{\CP}_1,{\CQ}_1]] + \half
  d_\CQ[[{\CP}_1,{\CP}_1]] &=
  [[{\CP},[[{\CP}_1,{\CQ}_1]]\,]] +
  \half[[{\CQ},[[{\CP}_1,{\CP}_1]]\,]] \\
  &=
  [[\,[[{\CP},{\CP}_1]],{\CQ}_1]]-[[{\CP}_1,[[{\CP},{\CQ}_1]]\,]]
  \\
  & \quad + \half [[\,[[{\CQ},{\CP}_1]],{\CP}_1]]-\half
  [[{\CP}_1,[[{\CQ},{\CP}_1]]\,]] \\
  &= [[{\CP}_1,[[{\CP},{\CQ}_1]]+[[{\CQ},{\CP}_1]]\,]] \\
  &= 0 .
\end{align*}

\subsection{Infinitesimal symmetries} \label{InfSymm}

Using Theorem \ref{Getzler}, we calculate the space
$H^0(\CL;d_{\CP},d_{\CQ})$ of infinitesimal symmetries of the
compatible pair of Hamiltonian operators $\CP$ and $\CQ$ of the
dispersionless KdV hierarchy. The Hamiltonian operator $\CP$ is
nondegenerate, so Theorem \ref{Getzler} applies to it.

Suppose that $Z\in Z^0(\CL\<\ell\>;d_\CP,d_\CQ)$, that is, $Z\in\CL^0$
and $d_\CP Z=d_\CQ Z=0$. Let $h\in\CA$ be the characteristic of $Z$.
If $\ell=0$, so that $h\in\CA_0$, we see that
\begin{equation*}
  d_\CQ \tint h(u) \theta \, dx = - \tint \bigl( uh'(u) + \half h(u)
  \bigr) \theta \theta_1 \, dx ,
\end{equation*}
which cannot be nonzero for smooth $h(u)$. Thus, we may assume that
$\ell>0$.

Theorem \ref{Getzler} implies that there exists a functional $\tint
g\,dx \in\CL^{-1}\<\ell\>$ such that
\begin{equation*}
  Z = d_\CP \tint g \, dx .
\end{equation*}
If $\ell=1$, we see that $g\in\CA_0$, and hence
\begin{equation*}
  Z = - \tint g'(u) \theta_1 \, dx = \tint u_1 g''(u) \theta \, dx .
\end{equation*}
If $\ell=2$, the density $g=u_1s(u)$ is a total derivative, and hence
$Z=0$.  Thus, assume that $\ell>2$.

Applying Theorem \ref{Getzler} once more to the equation
\begin{equation*}
  d_\CP \bigl( d_\CQ \tint g \, dx \bigr) = - d_\CQ d_\CP \tint g \,
  dx = 0 ,
\end{equation*}
we see that there is exists a functional $\tint
f\,dx\in\CL^{-1}\<\ell\>$ such that
\begin{equation} \label{fg}
  d_\CP \tint f \, dx = d_\CQ \tint g \, dx .
\end{equation}
We now prove that this is impossible.

The vectors $d_{\CP}\tint f \, dx = \tint \theta_1 \delta f \, dx$ and
$d_{\CQ}\tint g \, dx = \tint (u\theta_1+\half u_1\theta) \delta g \,
dx$ have characteristics $- \p \delta f$ and $- \p \delta(ug) + \half
u_1 \delta g$. Thus, \eqref{fg} may be written
\begin{equation*}
  \half u_1 \delta g = \p \bigl( \delta(ug) - \delta f \bigr) .
\end{equation*}
In particular, $u_1\delta g$ is a total derivative.
  
Suppose that $g\in\CA[n]$, $n>1$. We have
\begin{equation*}
  u_1 \delta g \equiv (-1)^n (1-n) u_1 u_{2n} \p^2_ng + \p( (-1)^n n u_1
  u_{2n-1}\p_n^2g ) \mod \CA[2n-2] .
\end{equation*}
This cannot be a total derivative unless $\p_n^2g=0$, in which case we
may replace $g$ by $g-\p(u_{n-1}\p_ng)\in\CA[n-1]$.
  
Arguing by downward induction, we may assume that $g\in\CA[1]$. In
other words, $g=u_1^{\ell-1}\gamma(u)$, where $\gamma$ is a smooth
function of $u$. We calculate that
\begin{equation}
  \label{deltaf1}
  \begin{aligned}
    \p \delta f &= \p \delta(ug) - \half u_1 \delta g \\
    &= - (\ell-2) \p \bigl( (\ell-1) u_2 u_1^{\ell-3} u \gamma(u) +
    u_1^{\ell-1} \bigl( u \gamma'(u) + \half \gamma(u) \bigr) \bigr) .
  \end{aligned}
\end{equation}
This implies that, up to a total derivative, $f$ has the form
$u_1^{\ell-1}\eta(u)$, and that
\begin{equation}
  \label{deltaf2}
  \delta f = - (\ell-2) \bigl( (\ell-1) u_2u_1^{\ell-3} \eta(u) +
  u_1^{\ell-1} \eta'(u) \bigr)
\end{equation}
Comparing \eqref{deltaf1} and \eqref{deltaf2}, we see that
\begin{equation*}
  (\ell-1) u_2 u_1^{\ell-3} u \gamma(u) + u_1^{\ell-1} \bigl( u
  \gamma'(u) + \half \gamma(u) \bigr)
  = (\ell-1) u_2u_1^{\ell-3} \eta(u) + u_1^{\ell-1} \eta'(u) ,
\end{equation*}
which has no nonzero solutions.

In summary, we have proven the following theorem.
\begin{theorem}
  The infinitesimal symmetries
  \begin{equation*}
    H^0(\CL;d_{\CP},d_{\CQ}) = \{ Z\in\CL^0 \mid d_{\CP}Z=d_{\CQ}Z=0 \}
  \end{equation*}
  of the bihamiltonian structure of the dispersionless KdV hierarchy 
  are the vector fields of the
  form $d_{\CP}\tint g(u)\,dx$, or equivalently, vector fields with
  characteristic of the form $u_1h(u)$.
\end{theorem}

\subsection{The bihamiltonian cohomology for structures of hydrodynamic type}

Let $\CP$ and $\CQ$ be compatible Hamiltonian operators of
hydrodynamic type. The bihamiltonian cohomology decomposes into
subspaces
\begin{equation*}
  H^k(\CL;d_{\CP},d_{\CQ}) = \bigoplus_{\ell\ge-k}
  H^k(\CL\<\ell\>;d_{\CP},d_{\CQ}) .
\end{equation*}

If $\CP$ is nondegenerate, any cohomology class in
$H^k(\CL\<\ell\>;d_\CP,d_\CQ)$, $\ell>-k$, has a representative of the
form $(0,\dots,0,c)$. The argument, which uses Theorem~\ref{Getzler},
is as follows.

Let $(c_0,\ldots,c_k)$ be a cocycle in
$Z^k(\CL\<\ell\>;d_{\CP},d_{\CQ})$, $\ell>-k$.  Since $d_{\CP}c_0=0$,
and $d_{\CP}$ is acyclic, we can write the first component $c_0$ of
the cocycle as a coboundary $d_{\CP}a_0$ for some
$a_0\in{\CL}^{k-1}\<\ell\>$. Subtracting
\begin{equation*}
  d(a_0,0,\dots,0) = (d_{\CP}a_0,d_{\CQ}a_0,0,\dots,0)
\end{equation*}
from $(c_0,\dots,c_k)$ gives a cocycle in the same cohomology class,
with vanishing first component:
\begin{equation*}
  (c_0,\ldots,c_k) \sim (0,c_1-d_{\CQ}a_0,c_2,\ldots,c_k) .
\end{equation*}
In turn, the equation $d_{\CP}(c_1-d_{\CQ}a_0)=0$ holds.  Iterating
the above procedure $k-1$ more times, one finally obtains a cocycle
$(0,\ldots,0,c)$ such that
\begin{equation*}
  (c_0,\ldots,c_k) \sim (0,\dots,0,c)
\end{equation*}
and $d_{\CP}c=d_{\CQ}c=0$. It follows that $c=d_{\CP}Y$ for some
$Y\in{\CL}^{k-1}\<\ell\>$. Since
\begin{equation*}
  d_{\CP}d_{\CQ}Y = - d_{\CQ}d_{\CP}Y = 0 ,
\end{equation*}
we may use the acyclicity of $d_{\CP}$ one last time to see the
existence of a cochain $X\in{\CL}^{k-1}\<\ell\>$ such that
$d_{\CQ}Y=d_{\CP}X$. In this way, we obtain the following.

\begin{proposition} \label{Hplus}
  Let $\CP$ and $\CQ$ be a pair of compatible Hamiltonian operators of
  hydrodynamic type such that $\CP$ is nondegenerate. If $\ell>-k$,
  the group $H^k(\CL\<\ell\>;d_{\CP},d_{\CQ})$ is isomorphic to
  \begin{equation*}
    \frac{ \bigl\{  d_{\CP}Y \mid \text{$Y\in\CL^{k-1}\<\ell\>$ and
        $d_{\CQ}Y=d_{\CP}X$ for some $X\in\CL^{k-1}\<\ell\>$} \bigr\}
    } { \bigl\{ d_{\CQ}T \mid \text{$T\in\CL^{k-1}\<\ell\>$ and
         $d_{\CP}T=0$} \bigr\} } .
   \end{equation*}
\end{proposition}

\begin{corollary} \label{CohGrps0and1}
  For the two Hamiltonian operators $\CP=\half\tint\theta\theta_1\,dx$
  and $\CQ=\half\tint u\theta\theta_1\,dx$ of the dispersionless KdV
  hierarchy, the group $H^k(\CL;d_{\CP},d_{\CQ})$ is isomorphic to
  \begin{equation*}
    \frac{ \bigl\{  d_{\CP}Y \mid \text{$Y\in\CL^{k-1}$ and
        $d_{\CQ}Y=d_{\CP}X$ for some $X\in\CL^{k-1}$} \bigr\} } {
      \bigl\{ d_{\CQ}T \mid \text{$T\in\CL^{k-1}$ and $d_{\CP}T=0$}
      \bigr\} } .
   \end{equation*}
\end{corollary}
\begin{proof}
  In light of Proposition~\ref{Hplus}, we must show that
  $H^k(\CL\<-k\>;d_{\CP},d_{\CQ})=0$ for $k\ge0$. In fact, the
  cohomology $H^k(\CL\<-k\>,d_\CP)$ is spanned by $\tint 1\,dx$ and
  $\tint u\,dx$ for $k=-1$, and by $\tint \theta\,dx$ for $k=0$. Thus,
  the above argument applies as well for $H^k(\CL\<-k\>,d_\CP)=0$, if
  $k>0$. 
  Finally, the case $k=0$ was discussed in Section~\ref{InfSymm}.  
\end{proof}

\subsection{Infinitesimal deformations} \label{QuasiTrivProof}

A bihamiltonian cohomology class $c\in H^\bullet(\CL;d_\CP,d_\CQ)$ is
quasi-trivial if its image in $H^\bullet(\CLhat;d_\CP,d_\CQ)$ is zero.
In other words, $c=(c_0,\dots,c_k)\in\CL^k\oplus\dots\oplus\CL^k$ is
quasi-trivial if there exists
$(b_0,\dots,b_{k-1})\in\CLhat^{k-1}\oplus\dots\oplus\CLhat^{k-1}$ such
that
\begin{equation*}
  d(b_0,\dots,b_{k-1}) = (c_0,\dots,c_k) .
\end{equation*}

We now state our main theorem. 
\begin{theorem} \label{mainthm}
  All infinitesimal bihamiltonian deformations of
  the bihamiltonian structure of the dispersionless
  KdV hierarchy of homogeneous degree greater than $1$ are
  quasi-trivial. In other words, the image of
  $H^1(\CL\<\ell\>;d_\CP,d_\CQ)$ in $H^1(\CLhat\<\ell\>;d_\CP,d_\CQ)$
  is zero if $\ell>0$.
\end{theorem}

Note that $H^1(\CL\<0\>;d_\CP,d_\CQ)$ is certainly not trivial: in
fact, it may be identified with the space of cocycles
\begin{equation*}
  \{ \bigl( 0,\tint s(u) \theta\theta_1 \, dx \bigr) \mid s(u) \in
  \CA_0 \}
\end{equation*}
modulo the one-dimensional space of coboundaries with basis
\begin{equation*}
  \bigl( 0 , \tint \theta\theta_1 \, dx \bigr) = \bigl( - d_\CP 2 \tint
  \theta \, dx , - d_\CQ 2 \tint \theta \, dx \bigr) .
\end{equation*}

As part of a more general study of infinitesimal bihamiltonian
deformations of ``semisimple'' bihamiltonian structures in any number
of dimensions, Liu and Zhang \cite{LZ1} proved that the kernel of the
map
\begin{equation*}
  H^1(\CL;d_\CP,d_\CQ) \to H^1(\CLhat;d_\CP,d_\CQ)
\end{equation*}
has the form
\begin{equation*}
  \{ (0,\tint s(u)\theta_1\theta_2\,dx) \mid s(u)\in\CA_0 \} .
\end{equation*}
The associated Hamiltonian operator may be calculated by applying the
normalization operator $\CN$:
\begin{equation*}
  \half \CN \tint s(u)\theta_1\theta_2\,dx = - \theta \bigl( s(u)\p^3 +
  \tfrac{3}{2} u_1s'(u)\p^2 + \half (u_2s'(u)+u_1^2s''(u)) \p \bigr)
  \theta .
\end{equation*}
In other words, the equivalence classes of infinitesimal deformations
of the bihamiltonian structure $\bigl( \p , u\p+\half u_1 \bigr)$ of
the dispersionless KdV hierarchy have the form
\begin{equation*} 
  \bigl( 0 , s(u)\p^3 + \tfrac{3}{2} u_1 s'(u) \p^2 + \half
  (u_2s'(u)+u_1^2s''(u)) \p \bigr) , \quad  s(u)\in\CA_0 .
\end{equation*}
When $s$ is constant, we recognize the infinitesimal deformation
associated to the bihamiltonian structure of the full KdV hierarchy.

\subsubsection{The cocycles}

We start the proof of Theorem \ref{mainthm} by studying the equation
\begin{equation} \label{FG}
  \theta\delta_\theta \bigl( d_\CP \tint f\theta\,dx - d_\CQ \tint
  g\theta\,dx \bigr) = 0 , \quad f,g \in \CA[n] .
\end{equation}

For $P=\half\theta\theta_1$ and $Q=\half u\theta\theta_1$, the differentials $d_\CP$ and $d_\CQ$ 
associated to the Hamiltonian operators $\CP=\tint P\,dx$ and $\CQ=\tint Q\,dx$ are given by the following formulas: if
$F\in\Lambda_\infty^\bullet$, or more generally, if
$F\in\hat\Lambda_\infty^\bullet$,
\begin{align*}
  d_\CP \tint F \, dx &= - \tint \delta_\theta P \, \delta_u F \, dx
  = - \tint \theta_1 \, \delta_u F \, dx \\
  &= -\sum_{k=0}^\infty \tint \theta_{k+1} \p_kF \, dx , \\
  d_\CQ \tint F\,dx &= - \tint \bigl( \delta_\theta Q \, \delta_u F +
  \delta_u Q \, \delta_\theta F \bigr) \, dx = - \tint \bigl( \bigl(
  u\theta_1+\half u_1\theta \bigr) \delta_u F +
  \half \theta \theta_1 \delta_\theta F \bigr) \, dx \\
  &= - \half \sum_{k=0}^\infty \tint \bigl( \bigl( \p^k(u\theta_1) +
  \p^{k+1}(u\theta) \bigr) \p_k F + \p^k(\theta\theta_1)
  \p_{\theta_k}F \bigr) \, dx .
\end{align*}
In particular, if $h\in\CAhat$, $d_\CP \tint h\,dx$ has characteristic
$\p\delta_uh$, and $d_\CQ \tint h\,dx$ has characteristic $\bigl(
u\p+\half u_1 \bigr) \delta_u h$.

\begin{proposition}\label{bivec}
  Let $f,g\in\CAhat[n]$. For $k\geq 0$, let
  \begin{equation*}
    F_k = \p_k f , \quad\text{and}\quad
    G_k = \half\sum_{\ell=0}^{n-k}\left[
      \tbinom{k+\ell}{\ell}+\tbinom{k+\ell+1}{\ell}\right]u_\ell\p_{k+\ell}
    g -\half\delta_{k,0}g.
  \end{equation*}
  Then
  \begin{align*}
    \theta\delta_\theta \bigl( d_\CP\tint f\theta\,dx \bigr) & =
    \sum_{k=0}^n \theta \theta_{k+1} \left( F_k + \sum_{j=k}^n (-1)^j
      \tbinom{j+1}{k+1} \p^{j-k} F_j \right) \intertext{and}
    \theta\delta_\theta \bigl( d_\CQ\tint g\theta\,dx \bigr) &=
    \sum_{k=0}^n \theta \theta_{k+1} \left( G_k + \sum_{j=k}^n (-1)^j
      \tbinom{j+1}{k+1} \p^{j-k} G_j \right) .
  \end{align*}
\end{proposition}
\begin{proof}
  We have
  \begin{align*}
    \theta\delta_\theta\left( d_\CP\tint f\theta\,dx \right) &=
    \theta\delta_\theta \left( \sum_{k=0}^n \theta\theta_{k+1} F_k
    \right) \\
    &= \theta \sum_{k=0}^n \bigl( \theta_{k+1}F_k
    + (-1)^k \p^{k+1}\left(\theta F_k\right) \bigr) \\
    &= \sum_{k=0}^n \theta \theta_{k+1} \left( F_k + \sum_{j=k}^n
      (-1)^j \tbinom{j+1}{k+1} \p^{j-k} F_j \right) .
  \end{align*}
  The formula for $\theta\delta_\theta \bigl( d_\CQ\tint g\,dx \bigr)$
  is derived similarly:
  \begin{align*}
    \theta\delta_\theta \bigl( d_\CQ\tint g\,dx \bigr) &= \half
    \theta\delta_\theta \left( \theta \sum_{k=0}^n \left(
        \p^k(u\theta_1) + \p^{k+1}(u\theta)\right) \p_k g -
      \theta\theta_1g \right) \\
    &= \half \theta\delta_\theta\left( \theta \sum_{k=0}^n
      \sum_{\ell=0}^k \left[
        \tbinom{k}{\ell}+\tbinom{k+1}{\ell}\right]
      \theta_{k-\ell+1} u_\ell\p_k g - \theta\theta_1 g \right) \\
    &= \half \theta\delta_\theta \left( \sum_{k=0}^n \theta
      \theta_{k+1} \sum_{\ell=0}^{n-k} \left[
        \tbinom{k+\ell}{\ell}+\tbinom{k+\ell+1}{\ell} \right]
      u_\ell\p_{k+\ell} g - \theta\theta_1g \right) \\
    &= \theta\delta_\theta \left( \sum_{k=0}^n \theta \theta_{k+1} G_k
    \right) \\
    &= \sum_{k=0}^n \theta \theta_{k+1} \left( G_k + \sum_{j=k}^n
      (-1)^j \tbinom{j+1}{k+1} \p^{j-k} G_j \right) .  \qedhere
  \end{align*}
\end{proof}

\subsubsection{The constraints}

For $j\geq 0$, let us define coefficients
\begin{align*}
  e_j &= F_j-G_j \\
  &= \p_j(f-ug) - \half \sum_{\ell=1}^{n-j} \left[
    \tbinom{j+\ell}{\ell}+\tbinom{j+\ell+1}{\ell} \right]
  u_\ell\p_{j+\ell} g + \tfrac{3}{2} \delta_{j,0} g ,
\end{align*}
and for $k\geq 0$ let
\begin{equation} \label{Sk}
  S_k = e_k + \sum_{j=0}^\infty (-1)^j \tbinom{j+1}{k+1} \p^{j-k}e_j .
\end{equation}
Assume that $f,g\in{\CAhat}[n]$. Then $e_j=0$ for $j>n$, and by
Proposition~\ref{bivec},
\begin{equation*}
  \theta\delta_\theta \left( d_\CP\tint f\theta\,dx - d_\CQ\tint
    g\theta\,dx \right) = \sum_{k=0}^n \theta\theta_{k+1} S_k .
\end{equation*}
That is, \eqref{FG} is equivalent to the system
\begin{equation*}
  {\CS} = \{S_k=0 \mid 0\leq k\leq n\} .
\end{equation*}

Let $n=2m>0$ be even. Define
\begin{equation}\label{E2i}
 E_\ell = \sum_{j=2\ell}^{m+\ell} (-1)^j \tbinom{2m-j}{m-\ell}
 \tbinom{j+1}{2\ell+1} \p^{j-2\ell} e_j .
\end{equation}
\begin{proposition}
  \label{theEs}
  We have
  \begin{equation*}
    S_k = \begin{cases}
      \displaystyle \sum_{\ell=0}^m
      \tfrac{\tbinom{2\ell+1}{k+1}}{\tbinom{2m-k-1}{m-\ell}} \,
      \p^{2\ell-k}E_\ell , & k<n , \\[10pt]
      2E_m , & k=n .
    \end{cases}
  \end{equation*}
  In particular, the subset $\{S_{2i}=0 \mid 0\leq i\leq m\}$ of the
  system $\CS$ of equations is equivalent to the system of equations
  $\CE=\{E_\ell=0 \mid 0\le \ell\le m\}$.
\end{proposition}
\begin{proof}
  The result is clear for $k=n$, since in this case, $S_n=2e_n$ and
  $E_m=e_n$. From now on, we assume that $k<n$.
  
  We wish to prove that
  \begin{align*}
    S_k = \sum_{\ell=0}^m
    \tfrac{\binom{2\ell+1}{k+1}}{\binom{n-k-1}{m-\ell}} \,
    \p^{2\ell-k} E_\ell &= \sum_{\ell=0}^m
    \tfrac{\binom{2\ell+1}{k+1}}{\binom{n-k-1}{m-\ell}} \,
    \sum_{j=2\ell}^{m+\ell} (-1)^j \tbinom{n-j}{m-\ell}
    \tbinom{j+1}{2\ell+1} \p^{j-k} e_j \\
    &= \sum_{j=k}^n (-1)^j \left( \sum_{\ell=0}^\infty
      \tfrac{\tbinom{n-j}{m-\ell}}{\binom{n-k-1}{m-\ell}}
      \tbinom{2\ell+1}{k+1} \tbinom{j+1}{2\ell+1} \right) \p^{j-k} e_j
    .
  \end{align*}
  (In this formula, it is understood that $\tbinom{r}{s}=0$ if $s<0$.)
  Comparing with the definition of $S_k$, we see that we are left to
  prove that for $k\le j\le n$,
  \begin{equation*}
    \sum_{\ell=0}^\infty
    \tfrac{\tbinom{n-j}{m-\ell}}{\binom{n-k-1}{m-\ell}}
    \tbinom{2\ell+1}{k+1} \tbinom{j+1}{2\ell+1} =
    \begin{cases} \tbinom{j+1}{k+1} , & j>k , \\
      1 + (-1)^k , & j=k .
    \end{cases}
  \end{equation*}
  We start with the case where $j>k$. Making the substitutions
  $\alpha=j-k-1$, $\beta=n-j$ and $\ell=p+j-m$, we have
  \begin{align*}
    \sum_{\ell=0}^\infty
    \frac{\tbinom{n-j}{m-\ell}}{\binom{n-k-1}{m-\ell}}
    \tbinom{2\ell+1}{k+1} \tbinom{j+1}{2\ell+1} &= \sum_{p=0}^\infty
    \frac{\tbinom{\beta}{\beta-p}}{\binom{\alpha+\beta}{\beta-p}}
    \tbinom{2p+2j-n+1}{k+1} \tbinom{j+1}{2p+2j-n+1} \\
    &= \frac{\beta!(j+1)!}{(\alpha+\beta)!(k+1)!}  \sum_{p=0}^\infty
    \frac{(\alpha+p)!}{p!(\alpha-\beta+2p+1)!(\beta-2p)!} \\
    &= \frac{\beta!(j+1)!}{(\alpha+\beta)!(k+1)!(\alpha+1)}
    \sum_{p=0}^\infty \tbinom{\alpha+p}{\alpha}
    \tbinom{\alpha+1}{\beta-2p} .
  \end{align*}
  By Lemma \ref{Binom}, this equals $\tbinom{j+1}{k+1}$.
  
  If $j=k$, then the only possibly nonzero term of the sum
  \begin{equation*}
    \sum_{\ell=0}^\infty
    \frac{\tbinom{n-j}{m-\ell}}{\binom{n-k-1}{m-\ell}}
    \tbinom{2\ell+1}{k+1} \tbinom{j+1}{2\ell+1}
  \end{equation*}
  is that with $k=2\ell$, in which case it equals $2$; if $k$ is odd,
  the sum vanishes.
\end{proof}

\begin{lemma} \label{Binom}
  \begin{equation*}
    \sum_{p=0}^{\infty} \tbinom{\alpha+1}{\beta-2p}
    \tbinom{\alpha+p}{p} = \tbinom{\alpha+\beta}{\beta}
  \end{equation*}
\end{lemma}
\begin{proof}
\ignore{
  Using the formula
  \begin{equation*}
    \tbinom{n+k}{k} = (-1)^k \tbinom{-n-1}{k} ,
  \end{equation*}
  we see that we must prove that
  \begin{equation*}
    \sum_{p=0}^{\infty} (-1)^p \tbinom{\alpha+1}{\beta-2p}
    \tbinom{-\alpha-1}{p} = (-1)^\beta \tbinom{-\alpha-1}{\beta} .
  \end{equation*}
  Multiply the left-hand side by $x^\beta$ and sum over $\beta$, and
  we obtain
  \begin{equation*}
    (1+x)^{\alpha+1} \sum_{p=0}^{\infty} (-x^2)^p 
    \tbinom{-\alpha-1}{p} = (1+x)^{\alpha+1} (1-x^2)^{-\alpha-1} =
    (1-x)^{-\alpha-1} .
  \end{equation*}
  Likewise, multiplying the right-hand side by $x^\beta$ and summing
  over $\beta$, we obtain $(1-x)^{-\alpha-1}$. Comparing coefficients,
  the lemma follows.
}
Using the formula
\begin{equation*}
    \tbinom{n+k}{k} = (-1)^k \tbinom{-n-1}{k} ,
  \end{equation*}
  we see that we must prove that
  \begin{equation*}
    \sum_{p=0}^{\infty} (-1)^p \tbinom{\alpha+1}{\beta-2p}
    \tbinom{-\alpha-1}{p} = (-1)^\beta \tbinom{-\alpha-1}{\beta} .
  \end{equation*}
Expanding both sides of the identity
\begin{equation*}
(1+x)^{\alpha+1}(1-x^2)^{-\alpha-1}=(1-x)^{-\alpha-1}
\end{equation*}
as power series of $x$ and extracting the coefficient of 
$x^\beta$, the lemma follows.
\end{proof}

\subsubsection{The induction}\label{induction}

The next result is the main part of the proof of
Theorem~\ref{mainthm}.
\begin{theorem}\label{mainpart}
  Let $n=2m>4$. Suppose $f$ and $g$ are characteristics in $\CAhat[n]$
  which satisfy the equation
  \begin{equation*}
    d_{\CP}\tint f\theta\,dx = d_{\CQ}\tint g\theta\,dx .
  \end{equation*}
  Then there are densities $a,b,c\in\CAhat[m]$ such that
  $\tilde{f}=f+\p\,\delta_ua+(u\p+\half u_1)\delta_ub$ and
  $\tilde{g}=g-\p\,\delta_ub+(u\p+\half u_1)\delta_uc$ lie in
  $\CAhat[n-2]$.
\end{theorem}

The point of this theorem is that if $f$ and $g$ satisfy
the equation $d_{\CP}\tint f\theta\,dx = d_{\CQ}\tint g\theta\,dx$,
then so do $\tilde{f}$ and $\tilde{g}$.

This theorem will be proven in a number of steps.

\paragraph*{Step 1}

Since
\begin{equation*}
  E_m = e_n
\end{equation*}
vanishes, we have $f-ug\in\CAhat[n-1]$ or, equivalently,
$e_n=\p_n(f-ug)=0$,

\paragraph*{Step 2}

We now argue by induction on $1\le t\le m$ that
\begin{equation*}
  \p_{n-t+1}\p_ng = 0 ,
\end{equation*}
or equivalently, that $g = u_ng_0 + g_1$, where $g_0\in{\CAhat}[n-t]$ and
$g_1\in\CAhat[n-1]$. Assume as induction hypothesis that $\p_{n-s+1}\p_n
g=0$ for $1\le s<t$. Since $e_j\in\CAhat[n]$, we see that
\begin{equation*}
  E_{m-t} = \sum_{j=n-2t}^{n-t} (-1)^j\tbinom{n-j}{t}
  \tbinom{j+1}{n-2t+1}\p^{j+2t-n}e_j \in \CAhat[n+t] .
\end{equation*}
We now argue as follows: the coefficient of $u_{n+t}$ in $E_{m-t}$
equals
\begin{align*}
  [u_{n+t}]E_{m-t} &= (-1)^{n-t} \tbinom{n-t+1}{n-2t+1}
  [u_{n+t}]\p^t e_{n-t} \\
  &= (-1)^{n-t} (n-t+\tfrac{3}{2}) \tbinom{n-t+1}{n-2t+1}
  u_1\p_{n-t+1}\p_n g .
\end{align*}
Since $(n-t+\tfrac{3}{2})\tbinom{n-t+1}{n-2t+1}\neq 0$, we see that
$\p_{n-t+1}\p_n g=0$.

\paragraph{Step 3}

At this point, we know that $g=u_ng_0+g_1$, where $g_0\in\CAhat[m]$ and
$g_1\in\CAhat[n-1]$. 

\begin{lemma} \label{evencob}
  If $h\in\CAhat[m]$, $m>0$, then
  \begin{equation*}
    \p\,\delta_u(uh) - (u\p+\half u_1)\delta_uh \equiv (-1)^m
    (m+\half) u_nu_1\p_m^2h \mod \CAhat[n-1] .
  \end{equation*}
\end{lemma}
\begin{proof}
  We have
  \begin{align*}
    \delta_u(uh) &\equiv (-\p)^m u\p_mh \equiv u \delta_uh + (-1)^m m
    u_{n-1} u_1 \p_m^2h \mod \CAhat[n-2] , \\
    \delta_uh &\equiv (-\p)^m \p_mh \equiv (-1)^m u_n \p_m^2h \mod
    \CAhat[n-1] .
  \end{align*}
  From these equations, the lemma follows easily.
\end{proof}

Let $h\in\CAhat[m]$ be a solution of the equation
\begin{equation*}
  g_0 = (-1)^m (m+\half) u_1\p_m^2h .
\end{equation*}
By Lemma~\ref{evencob}, $f+\p\,\delta_u(u^2h)-(u\p+\half u_1)
\delta_u(uh)$ and $g+\p\,\delta_u(uh)-(u\p+\half u_1)\delta_uh$ are in
$\CAhat[n-1]$.

\paragraph{Step 4}

In Step 3, we have shown that we may reduce to the case that $f$,
$g\in\CAhat[n-1]$. We now show that $\p_{n-1}^2(f-ug)=0$. We have
\begin{equation*}
  E_{m-1} = 2 e_{n-2} - n\p e_{n-1} = 0 .
\end{equation*}
Since $e_{n-2}\in\CAhat[n-1]$ and $e_{n-1}=\p_{n-1}(f-ug)$, we see
that
\begin{align*}
  [u_n]E_{m-1} & = -n[u_n]\p e_{n-1} \\
  &= -n\p^2_{n-1}(f-ug)=0.
\end{align*}
In particular, $f-(ug+u_{n-1}\p_{n-1}(f-ug))\in\CAhat[n-2]$.

\paragraph{Step 5}

We now argue by induction on $2\le t\le m$ that 
\begin{equation*}
  \p_{n-t}\p_{n-1}(f-ug) = 0 ,
\end{equation*}
or equivalently, that $e_{n-1}=\p_{n-1}(f-ug)\in\CAhat[n-t-1]$. Assume
as induction hypothesis that $\p_{n-s}e_{n-1}=0$ for $1\le s < t$, so
that
\begin{equation*}
  \p_{n-t}e_{n-1}=[u_{2n-t-1}]\p^{n-1}e_{n-1}.
\end{equation*}
\begin{lemma}\label{theTs}
  If for $2\le t\le m$, and $1\le s < t$, $\p_{n-s}e_{n-1}=0$, then
  \begin{equation*}
    [u_{2n-t-1}]\p^{n-1}e_{n-1}=0.
  \end{equation*}
\end{lemma}
\begin{proof}
  For $2\leq t\leq m$ and $1\leq k\leq t$, let
  \begin{equation*}
    T_{t,k} = [u_{2n-t-1}] \sum_{j=n-t}^{n-1} \tbinom{j+1}{k-1}
    (-\p)^je_j .
  \end{equation*}
  Then 
  \begin{equation*}
    \sum_{k=1}^{n}(-1)^{k+t} \tbinom{n-k}{n-t} T_{t,k} = -
    [u_{2n-t-1}] \p^{n-1}e_{n-1}
  \end{equation*}
  since
  \begin{equation*}
    \sum_{k=1}^n(-1)^{j+k+t} \tbinom{n-k}{n-t} \tbinom{j+1}{k-1} =
    \begin{cases}
      0, & n-t< j<n-1, \\
      -1, & j=n-1.
    \end{cases}
  \end{equation*}
  From the definition~ \eqref{Sk} of the $S_k$, we see that
  \begin{equation*}
    T_{t,k} =
    \begin{cases}
      2 [u_{2n-t-1}] \p^{k-2} S_{k-2} , & 1< k\leq t , \\
      [u_{2n-t-1}] \sum_{j=0}^{n-2} (-\p)^j S_j , & k=0.
    \end{cases}
  \end{equation*}
  Since the functions $S_k$ vanish, $T_{t,k}$ also vanish. It follows
  that $[u_{2n-t-1}]\p^{n-1}e_{n-1}=0$ for $2\leq t\leq m$.
\end{proof}

\paragraph{Step 6}

In Step 5, we showed that after the redefinition of Step 3, we have
$f=ug+u_{n-1}e_{n-1}+ f_0$, where $f_0\in\CAhat[n-2]$ and
$e_{n-1}\in\CAhat[m-1]$. The next lemma shows that after a further
redefinition of $f$, we may assume $e_{n-1}=0$, that is,
$f-ug\in\CAhat[n-2]$.

\begin{lemma}\label{oddcob}
  If $h\in\CAhat[m-1]$, $m>0$, then
  \begin{equation*}
    \delta_uh \equiv (-1)^{m-1} u_{n-2}\p_{m-1}^2h \mod \CAhat[n-3] .
  \end{equation*}
\end{lemma}
\begin{proof}
  We have
  \begin{equation*}
    \delta_uh \equiv (-\p)^{m-1} \p_{m-1}h \mod \CAhat[n-3] ,
  \end{equation*}
  and the lemma follows.
\end{proof}

Let $h\in\CAhat[m-1]$ be a solution of the equation
\begin{equation*}
  e_{n-1} = (-1)^m \p_{m-1}^2h .
\end{equation*}
Replacing $f$ by $f+\p\,\delta_uh$, we see that $f-ug\in\CAhat[n-2]$.

\paragraph{Step 7} 

We now argue by induction on $1\le t< m$ that
\begin{equation*}
  \p_{n-t}\p_{n-1}g=0,
\end{equation*}
or equivalently, that $g=u_{n-1}g_0+g_1$, where $g_0\in\CAhat[n-t-1]$ and
$g_1\in\CAhat[n-2]$.  Assume as induction hypothesis that
$\p_{n-s}\p_{n-1}g=0$ for $1\le s < t$. Since $e_j\in\CAhat[n-1]$, we see
that
\begin{equation*}
  E_{m-t-1}= \sum_{j=n-2t-2}^{n-t-1} (-1)^j \tbinom{n-j}{t+1}
  \tbinom{j+1}{n-2t-1} \p^{j+2t-n+2}e_j \in \CAhat[n+t] .
\end{equation*}
The coefficient of $u_{n+t}$ in $E_{m-t-1}$ equals
\begin{align*}
  [u_{n+t}]E_{m-t-1} &=
  (-1)^{n-t-1}\tbinom{n-t}{n-2t-1}[u_{n+t}]\p^{t+1}e_{n-t-1} \\
  &= (-1)^{n-t-1}(n-t+\half)\tbinom{n-t}{n-2t-1}u_1\p_{n-t}\p_{n-1}g.
\end{align*}
Since $(n-t+\half)\tbinom{n-t}{n-2t-1} \neq 0$, we see that
$\p_{n-t}\p_{n-1}g=0$ for $1\le t <m$.

\paragraph{Step 8}

In Step 7, we showed that $g=u_{n-1}g_0+g_1$, where $g_0\in\CAhat[m]$
and $g_1\in\CAhat[n-2]$. We now show that $g_0$ is actually in
$\hat\CA[m-1]$, that is, $\p_m\p_{n-1}g=0$.
\begin{lemma}
  If $e_{n-1}=0$, then $[u_{3m-2}]\p^{m-1}e_{m-1}=0$.
\end{lemma}
\begin{proof}
  For $1\leq k\le m$, let
  \begin{equation*}
    U_{m,k} = [u_{3m-2}] \sum_{j=m-1}^{n-1} \tbinom{j+1}{k-1}
    (-\p)^je_j .
  \end{equation*}
  Then 
  \begin{equation*}
    \sum_{k=1}^n (-1)^{k+m} \tbinom{n-k}{m-k} U_{m,k} =
    [u_{3m-2}]\p^{m-1}e_{m-1}
  \end{equation*}
  since $e_{n-1}=0$, and
  \begin{equation*}
    \sum_{k=1}^n (-1)^{j+k+m} \tbinom{n-k}{m-k} \tbinom{j+1}{k-1} =
    \begin{cases}
      1 , & j=m-1, \\
      0 , & m\le j < n-1 .
    \end{cases}
  \end{equation*}
  From the definition of $S_k$ we see that 
  \begin{equation*}
    U_{m,k} =
    \begin{cases}
      2[u_{3m-2}]\p^{k-2}S_{k-2} , & 0< k\leq m , \\
      [u_{3m-2}] \sum_{j=0}^{n-2} (-\p)^j S_j , & k=0.
    \end{cases}
  \end{equation*}
  Since the functions $S_k$ vanish, the lemma follows.
\end{proof}

To show that $\p_m\p_{n-1}g=0$, we now argue as follows.  By
definition,
\begin{equation*}
  e_{m-1} = \p_{m-1}(f-ug) - \half\sum_{\ell=1}^m \left[
    \tbinom{m+\ell-1}{\ell}+\tbinom{m+\ell}{\ell} \right]
  u_\ell\p_{m+\ell-1} g .
\end{equation*}
Since $f-ug\in\CAhat[n-2]$ and $g-u_ng_0\in\CAhat[n-2]$, we see that
the only contribution to the coefficient of $u_{3m-2}$ in
$\p^{m-1}e_{m-1}$ comes from the term in $e_{m-1}$ with $\ell=1$, and
that
\begin{equation*}
  [u_{3m-2}]\p^{m-1}e_{m-1} = -(m+\half) u_1 \p_m\p_{n-1}g .
\end{equation*}
Since $(m+\half)\neq 0$, we see that $\p_m\p_{n-1}g=0$.

\paragraph{Step 9}

We have shown that after redefinitions of $f$ and $g$, we have
$g=u_{n-1}g_0+g_1$, where $g_0\in\CAhat[m-1]$ and $g_1\in\CAhat[n-2]$,
and that $f-ug\in\CAhat[n-2]$. Let $h\in\CAhat[m-1]$ be a solution of
the equation
\begin{equation*}
  g_0 = (-1)^m \p_{m-1}^2h .
\end{equation*}
By Lemma~\ref{oddcob}, $f-(u\p+\half u_1)\delta_uh$ and
$g+\p\,\delta_uh$ are in $\CAhat[n-2]$.

This completes the proof of Theorem~\ref{mainpart}.

\subsubsection{The proof of Theorem \ref{mainthm}} \label{two}

Let $c=(c_0,c_1)$ be a cohomology class in $H^1(\CL;d_\CP,d_\CQ)$. By
Corollary \ref{CohGrps0and1}, there exists $n\ge0$ and $f$ and $g$ in
$\CA[n]$ such that $(c_0,c_1)$ is cohomologous to $(0,d_\CP\tint
g\theta\,dx)$ and the equation
\begin{equation*}
  d_{\CP}\tint f\theta\,dx = d_{\CQ}\tint g\theta\,dx
\end{equation*}
holds. We may assume that there exists an integer $\ell$ such that $c$
is homogeneous of degree $\ell+1$, that is, that $c\in
H^1(\CL\<\ell\>;d_\CP,d_\CQ)$. Then $f$ and $g$ may be taken to be
homogeneous of degree $\ell$; since they are polynomial in the jet
variables $\{u_1,\dots,u_n\}$, we conclude that $n$ is no larger than
$\ell$.

If $\ell>2$, we may redefine $f$ and $g$ so that they lie in
$\CAhat[2]$. To see this, we use a downward induction based on Theorem
\ref{mainpart} to redefine $f$ and $g$ so that they lie in
$\CAhat[4]$. All the steps up until Step 9 in
Subsection~\ref{induction} remain valid for $n=4$, showing that after
a further redefinition, we may assume that $f=ug+f_1$ and that
$g=u_3u_1^{\ell-3}s(u)+g_1$, where $f_1,g_1\in\CAhat[2]$. Since
$\ell>2$, the argument in Step 9 may still be used, and in this way,
we may redefine $f$ and $g$ so that they lie in $\CAhat[2]$.

It is easily checked that Steps 1 and 2 in Subsection~\ref{induction}
apply, and since $\ell>2$, Step 3 applies as well. Thus, we are
reduced to the case where $f$ and $g$ lie in $\CAhat[1]$, and hence
\begin{equation*}
  f=u_1^\ell s(u) \quad\text{and}\quad g=u_1^\ell t(u) ,
\end{equation*}
where $s,t\in\CA_0$. The equation $S_0=0$ may be rewritten as
\begin{equation*}
  0 = e_0 - \p e_1 = \bigl( \p_0(f-ug) - \tfrac{3}{2} u_1\p_1g +
  \tfrac{3}{2} g \bigr) - \p \, \p_1(f-ug) .
\end{equation*}
Taking the coefficient of $u_2$, we see that
\begin{align*}
  0 &= [u_2] ( e_0 - \p e_1 ) = - [u_2] \p \, \p_1(f-ug) \\
  &= - \ell(\ell-1)u_1^{\ell-2}(s(u)-ut(u)) .
\end{align*}
Thus, $f=ug$, and hence
\begin{equation*}
  e_0 - \p e_1 = \tfrac{3}{2} (1-\ell) g .
\end{equation*}
Therefore $g$ vanishes, and hence so does the
cohomology class $(0,d_\CP\tint g\theta\,dx)$.

We now turn to the case $\ell=2$. By the vanishing of
$e_2=\p_2(f-ug)$, we see that
\begin{equation*}
  f = ug + u_1^2 p(u) \quad\text{and}\quad g = u_2 s(u) + u_1^2 t(u) ,
\end{equation*}
with $s,t,p\in\CA_0$. In this case, the equation $S_0=0$ becomes
\begin{align*}
  0 &= e_0 - \p e_1 = \bigl( \p_0(f-ug) - \tfrac{3}{2} u_1\p_1g - 2
  u_2 \p_2 g + \tfrac{3}{2} g \bigr) - \p \bigl( \p_1(f-ug) -
  \tfrac{5}{2} u_1\p_2g \bigr) \\
  &= 2 u_2 ( s(u) - p(u) ) + u_1^2 \bigl( \tfrac{5}{2} s'(u) - p'(u) -
  \tfrac{3}{2} t(u) \bigr) .
\end{align*}
It follows that $s(u)=p(u)$ and that $t(u)=p'(u)$, and hence that
$g=\p(u_1p(u))$. We calculate that
\begin{equation*}
  d_\CP \tint g\theta \, dx = - d_\CP \tint u_1 p(u) \theta_1 \, dx =
  - \tint p(u) \theta_1\theta_2 \, dx .
\end{equation*}
The following lemma shows that the cocycle $(0,\tint p(u)
\theta_1\theta_2 \, dx)$ is quasi-trivial.
\begin{lemma}
  \begin{equation*}
    d_{\CP} d_{\CQ} \int \frac{u_2}{u_1} \, h(u) \, dx =
    \frac{3}{2} \int h'(u) \,\theta_1 \theta_2 \, dx 
  \end{equation*}
\end{lemma}
\begin{proof}
  We calculate that
  \begin{equation*}
    d_{\CQ} \int \frac{u_2}{u_1} \, h(u) \, dx = - \int \biggl[
    \frac{3}{2} h(u) \theta_2 + uh'(u) \Bigl( \frac{u_2}{u_1} \,
    \theta_1 - \theta_2 \Bigr) \biggr] \, dx .
  \end{equation*}
  Applying $d_{\CP}$, the lemma follows.
\end{proof}

The case $\ell=1$ is uninteresting: since $g=u_1s'(u)$, $s(u)\in\CA_0$,
we have
\begin{equation*}
  d_\CP \tint g\theta \, dx = - d_\CP \tint s(u) \theta_1 \, dx = 0 .
\end{equation*}

\pagebreak

\end{document}